\documentclass[11pt,leqno]{amsart}
\usepackage{amsmath}
\usepackage{amsthm}
\usepackage{amssymb}
\theoremstyle{definition}

\theoremstyle{plain}

\begin{document}
\setcounter{page}{1}
\begin{flushleft}
\scriptsize{{Inequality Theory and Applications\/} {\bf 6}\,(2007), %\hfill\hfil ti dkdkdkd
\\
Pages ???--???}
\end{flushleft}
\vspace{16mm}

\begin{center}
{\normalsize\bf On $H_3(1)$ Hankel determinant for some classes of univalent functions } \\[12mm]
     {\sc K. O. BABALOLA$^{1,2}$} \\ [8mm]
%     {\small  Submitted March 19, 2005 }\\[5mm]

\begin{minipage}{123mm}
{\small {\sc Abstract.}
   Focus in this paper is on the Hankel determinant, $H_3(1)$, for the well-known classes of bounded-turning, starlike and convex functions in the open unit disk $E=\{z\in \mathbb{C}\colon|z|<1\}$. The results obtained complete the series of research works in the search for sharp upper bounds on $H_3(1)$ for each of these classes.}
\end{minipage}
\end{center}

 \renewcommand{\thefootnote}{}
 %\footnotetext{Received September 22, 2005.\qquad $^\ast$\,Corresponding author.}
 \footnotetext{2000 {\it Mathematics Subject Classification.}
            30C45; 30C50.}
 \footnotetext{{\it Key words and phrases.} Hankel determinant, functions of bounded turning, starlike and convex univalent functions.}
 %\footnotetext{This paper was supported by .....}
 \footnotetext{$^1$\,Current Address: Centre for Advanced Studies in Mathematics,
Lahore University of Management Sciences, Lahore, Pakistan. E-mail:
kobabalola\symbol{64}lums.edu.pk}
 \footnotetext{$^2$\,Permanent Address: Department of Mathematics, University of Ilorin,
 Ilorin, Nigeria. E-mail:
kobabalola\symbol{64}gmail.com; babalola.ko\symbol{64}unilorin.edu.ng}

\def\iff{if and only if }
\def\S{Smarandache }
\def\pi{positive implicative }
\newcommand{\norm}[1]{\left\Vert#1\right\Vert}
% ----------------------------------------------------------------

%\section{}
\vskip 12mm
{\bf 1. Introduction }
\medskip

Let $A$ be the class of functions
$$f(z)=z+a_2z^2+\cdots\eqno{(1.1)}$$
which are analytic in $E$. A function $f\in A$ is said to be of bounded turning, starlike and convex respectively if and only if, for $z\in E$, Re $f'(z)>0$, Re $zf'(z)/f(z)>0$ and Re $(1+zf''(z)/f'(z))>0$. By usual notations we denote these classes of functions respectively by $R$, $S^\ast$ and $C$. Let $n\geq 0$ and $q\geq 1$, the $q$-th Hankel determinant is defined as:
\[
H_q(n)=
\begin{vmatrix}
a_n & a_{n+1} & \cdots & a_{n+q-1}\\
a_{n+1} & \cdots & \cdots & \vdots\\
\vdots & \vdots & \vdots & \vdots\\
a_{n+q-1} & \cdots & \cdots & a_{n+2(q-1)}
\end{vmatrix}
\]
(see \cite{JW} for example). The determinant has been investigated by several authors with the subject of inquiry ranging from rate of growth of $H_q(n)$ as $n\rightarrow\infty$ \cite{JW,KI} to the determination of precise bounds on $H_q(n)$ for specific $q$ and $n$ for some favored classes of functions \cite{AJ,AS,AK}. In particlar, sharp upper bounds on $H_2(2)$ were obtained by the authors of articles \cite{AJ,AS,AK} for various classes of functions. In the present investigation, our focus is on the Hankel determinant, $H_3(1)$, for the well-known classes of bounded-turning, starlike and convex functions in $E$.\vskip 2mm

By definition, $H_3(1)$ is given by
$$H_3(1)=
\begin{vmatrix}
a_1 & a_2 & a_3\\
a_2 & a_3 & a_4\\
a_3 & a_4 & a_5
\end{vmatrix}.$$
For $f\in A$, $a_1=1$ so that
$$H_3(1)=a_3(a_2a_4-a_3^2)-a_4(a_4-a_2a_3)+a_5(a_3-a_2^2)$$
and by triangle inequality, we have
$$|H_3(1)|\leq|a_3||a_2a_4-a_3^2|+|a_4||a_2a_3-a_4|+|a_5||a_3-a_2^2|.\eqno{(1.2)}$$
\vskip 2mm

Incidentally, all of the functionals on the right side of the inequality (1.2) have known (and sharp) upper bounds in the classes of functions which are of interest in this paper, except $|a_2a_3-a_4|$. The last one is the well-known Fekete-Szego functional. For $R$, sharp bound $2/3$ was reported in \cite{BO} (with $R$ corresponding to $n=\alpha=1$, $\beta=0$ in the classes $T_n^\alpha(\beta)$ studied there) while for $S^\ast$ and $C$, sharp bounds $1$ and $1/3$ respectively were given in \cite{FR}. Janteng et-al \cite{AJ,AS} obtained for the functional $|H_2(2)|\equiv|a_2a_4-a_3^2|$ sharp bounds $4/9$, $1$ and $1/8$ repectively for $R$, $S^\ast$ and $C$. Furthermore, it is known that for $k=2,3,\cdots$, $|a_k|\leq 2/k$, $|a_k|\leq k$ and $|a_k|\leq 1$ also respectively for $R$, $S^\ast$ and $C$ (see \cite{PL,TH}). Thus finding the best possible bounds on $|a_2a_3-a_4|$ for each of the classes and using those known inequalities, then the sharp upper bounds on $H_3(1)$ follow as simple corollaries.\vskip 2mm

Our investigation follows a method of classical analysis devised by Libera and Zlotkiewicz \cite{EJ,RJ}. The same has been employed by many authors in similar works (see also \cite{AJ,AS,AK}). In the next section we state the necessary lemmas while in Section 3 we present our main results.
\bigskip

{\bf 2. Preliminary Lemmas}
\medskip

Let $P$ denote the class of functions $p(z)=1+c_1z+c_2z^2+\cdots$ which are regular in $E$ and satisfy Re $p(z)>0$, $z\in E$. To prove the main results in the next section we shall require the following two lemmas.\vskip 2mm

{\bf Lemma 2.1.} (\cite{PL}) {\em
Let $p\in P$, then $|c_k|\leq 2$, $k=1,2,\cdots$, and the inequality is sharp.}\vskip 2mm

{\bf Lemma 2.2.} (\cite{EJ,RJ}) {\em Let $p\in P$, then 
$$2c_2=c_1^2+x(4-c_1^2)\eqno{(2.1)}$$
and
$$4c_3=c_1^3+2xc_1(4-c_1^2)-x^2c_1(4-c_1^2)+2z(1-|x|^2)(4-c_1^2)\eqno{(2.2)}$$
for some $x$, $z$ such that $|x|\leq 1$ and $|z|\leq 1$.} 
\bigskip

{\bf 3. Main Results}
\medskip

{\bf Theorem 3.1.} {\em Let $f\in R$. Then
$$|a_2a_3-a_4|\leq\frac{1}{2}.$$
The inequality is sharp. Equality is attained by
$$f(z)=\int_0^z\frac{1+t^3}{1-t^3}dt.$$}\vskip 2mm
\begin{proof}
Let $f\in R$. Then there exists a $p\in P$ such that $f'(z)=p(z)$, wherefrom equating coefficients we find that $2a_2=c_1$, $3a_3=c_2$ and $4a_4=c_3$. Thus we have
$$|a_2a_3-a_4|=\left|\frac{c_1c_2}{6}-\frac{c_3}{4}\right|.\eqno{(3.1)}$$
Substituting for $c_2$ and $c_3$ using Lemma 2, we obtain
$$|a_2a_3-a_4|=\left|\frac{c_1^3}{48}-\frac{c_1(4-c_1^2)x}{24}+\frac{c_1(4-c_1^2)x^2}{16}-\frac{(4-c_1^2)(1-|x|^2)z}{8}\right|.\eqno{(3.2)}$$
By Lemma 1, $|c_1|\leq 2$. Then letting $c_1=c$, we may assume without restriction that $c\in [-2,0]$. Thus applying the triangle inequality on (3.2), with $\rho=|x|$, we obtain
$$\aligned |a_2a_3-a_4|
&\leq\frac{c^3}{48}+\frac{(4-c^2)}{8}+\frac{c(4-c^2)\rho}{24}+\frac{(c-2)(4-c^2)\rho^2}{16}\\
&=F(\rho).
\endaligned$$
Now we have
$$F'(\rho)=\frac{c(4-c^2)}{24}+\frac{(c-2)(4-c^2)\rho}{8}<0.$$
Hence $F(\rho)$ is a decreasing function of $\rho$ on the closed interval $[0,1]$, so that $F(\rho)\leq F(0)$. That is
$$\aligned F(\rho)
&\leq\frac{c^3}{48}+\frac{4-c^2}{8}\\
&=G(c).
\endaligned$$
Obviously $G(c)$ is increasing on $[-2,0]$. Hence we have $G(c)\leq G(0)=1/2$.\vskip 2mm

By setting $c_1=c=0$ and selecting $x=0$ and $z=1$ in (2.1) and (2.2) we find that $c_2=0$ and $c_3=2$. Thus equality is attained by $f(z)$ defined in theorem and the proof is complete.
\end{proof}\vskip 2mm

Let $f\in R$. Then using the above result in (1.2) together with the known inequalities $|a_3-a_2^2|\leq 2/3$ \cite{BO}, $|a_2a_4-a_3^2|\leq 4/9$ \cite{AJ} and $|a_k|\leq 2/k$, $k=2,3,\cdots$ \cite{TH}, we have the sharp inequality:\vskip 2mm

{\bf Corollary 3.2.} {\em Let $f\in R$. Then
$$|H_3(1)|\leq\frac{993}{1620}.$$}\vskip 2mm

{\bf Theorem 3.3.} {\em Let $f\in S^\ast$. Then
$$|a_2a_3-a_4|\leq 2.$$
The inequality is sharp. Equality is attained by the Koebe function $k(z)=z/(1-z)^2$.}\vskip 2mm
\begin{proof}
Let $f\in S^\ast$. Then there exists a $p\in P$ such that $zf'(z)=f(z)p(z)$. Equating coefficients we find that $a_2=c_1$, $2a_3=c_2+c_1^2$ and $6a_4=2c_3+3c_1c_2+c_1^3$. Thus we have
$$|a_2a_3-a_4|=\frac{1}{3}|c_1^3-c_3|.\eqno{(3.3)}$$
Substituting for $c_3$ from Lemma 2, we obtain
$$|a_2a_3-a_4|=\frac{1}{12}|3c_1^3-2c_1(4-c_1^2)x+c_1(4-c_1^2)x^2-2(4-c_1^2)(1-|x|^2)z|.\eqno{(3.4)}$$
Since $|c_1|\leq 2$ by Lemma 1, let $c_1=c$ and assume without restriction that $c\in [0,2]$. Applying the triangle inequality on (3.4), with $\rho=|x|$, we obtain
$$\aligned |a_2a_3-a_4|
&\leq\frac{1}{12}[3c^3+2(4-c^2)+2c(4-c^2)\rho+(c-2)(4-c^2)\rho^2]\\
&=F(\rho).
\endaligned$$
Differentiating $F(\rho)$, we have
$$F'(\rho)=\frac{1}{12}[2c(4-c^2)+2(c-2)(4-c^2)]>0.$$
This implies that $F(\rho)$ is an increasing function of $\rho$ on $[0,1]$ if $c\in[1,2]$. In this case $F(\rho)\leq F(1)=c\leq 2$ for all $\rho\in[0,1]$. It follows therefore that $F(\rho)\leq 2$. On the other hand suppose $c\in[0,1)$, then $F(\rho)$ is decreasing on $[0,1]$ so that $F(\rho)\leq F(0)$. That is
$$\aligned F(\rho)
&\leq\frac{3c^3-2c^2+8}{12}\\
&=G(c).\endaligned$$
Hence we have $G(c)\leq G(0)=2/3$, $c\in[0,1)$. This is less than 2, which is the case when $c\in[1,2]$. Thus the maximum of the functional $|a_2a_3-a_4|$ corresponds to $\rho=1$ and $c=2$.\vskip 2mm

If $c_1=c=2$ in (2.1) and (2.2), then we have $c_2=c_3=2$. Using these in (3.3) we see that equality is attained which shows that our result is sharp. Furthermore, it is easily seen that the extremal function in this case is the well known Koebe function $k(z)=z/(1-z)^2$.
\end{proof}\vskip 2mm

For $f\in S^\ast$, using the known inequalities $|a_k|\leq k$, $k=2,3,\cdots$ \cite{PL}, $|a_2a_4-a_3^2|\leq 1$ \cite{AS} and $|a_3-a_2^2|\leq 1$ \cite {FR} together with Theorem 2 we have the next corollary.\vskip 2mm

{\bf Corollary 3.4.} {\em Let $f\in S^\ast$. Then
$$|H_3(1)|\leq 16.$$
The inequality is sharp. Equality is attained by a rotation, $k_1(z)=z/(1+z)^2$, of the Koebe function.}\vskip 2mm

{\bf Theorem 3.5.} {\em Let $f\in C$. Then
$$|a_2a_3-a_4|\leq\frac{1}{6}.$$ The inequality is sharp. Equality is attained by
$$f(z)=\int_0^z\left\{s.\exp\left(\int_0^s\frac{2t^3}{1-t^3}dt\right)\right\}ds.$$
}\vskip 2mm
\begin{proof}
For $f\in C$ given by (1.1), there exists a $p\in P$ such that $(zf'(z))'=f'(z)p(z)$. Then equating coefficients we find that $2a_2=c_1$, $6a_3=c_2+c_1^2$ and $24a_4=2c_3+3c_1c_2+c_1^3$. Thus we have
$$|a_2a_3-a_4|=\frac{1}{24}|c_1^3-c_1c_2-2c_3|.\eqno{(3.5)}$$
Substituting for $c_2$ and $c_3$ using Lemma 2, we obtain
$$|a_2a_3-a_4|=\frac{1}{48}|-3c_1(4-c_1^2)x+c_1(4-c_1^2)x^2-2(4-c_1^2)(1-|x|^2)z|.\eqno{(3.6)}$$
With $|c_1|\leq 2$ from Lemma 1, we let $c_1=c$ and assume also without restriction that $c\in [-2,0]$. Thus applying the triangle inequality on (3.6), with $\rho=|x|$, we obtain
$$\aligned |a_2a_3-a_4|
&\leq\frac{(4-c^2)}{24}+\frac{c(4-c^2)\rho}{16}+\frac{(c-2)(4-c^2)\rho^2}{48}\\
&=F(\rho).
\endaligned$$
Differentiating $F(\rho)$, we get
$$F'(\rho)=\frac{c(4-c^2)}{16}+\frac{(c-2)(4-c^2)\rho}{24}<0.$$
Thus $F(\rho)$ is a decreasing function of $\rho$ on $[0,1]$, so that $F(\rho)\leq F(0)$. That is
$$\aligned F(\rho)
&\leq\frac{4-c^2}{24}\\
&=G(c),
\endaligned$$
which is increasing on $[-2,0]$. Hence $G(c)\leq G(0)=1/6$. Thus the maximum of the functional $|a_2a_3-a_4|$ corresponds to $c=0$ and $\rho=0$, which is $1/6$.\vskip 2mm

If we set $c_1=c=0$ and selecting $x=0$ and $z=1$ in (2.1) and (2.2) we find that $c_2=0$ and $c_3=2$, and equality is attained by $f(z)$ defined in theorem. This completes the proof.
\end{proof}\vskip 2mm

Finally for $f\in C$ if we use the known inequalities $|a_k|\leq 1$, $k=2,3,\cdots$ \cite{PL}, $|a_2a_4-a_3^2|\leq 1/8$ \cite{AS} and $|a_3-a_2^2|\leq 1/3$ \cite {FR} together with the last result, we obtain the following sharp inequality:\vskip 2mm

{\bf Corollary 3.6.} {\em Let $f\in C$. Then
$$|H_3(1)|\leq\frac{15}{24}.$$}

\medskip
{\it Acknowledgements.} This work was carried out at the Centre for Advanced Studies in Mathematics, CASM, Lahore University of Management Sciences, Lahore, Pakistan during the author's postdoctoral fellowship at the Centre. The author is indebted to all staff of CASM for their hospitality, most especially Prof. Ismat Beg. 
\bigskip

 \bibliographystyle{amsplain}

\end{document}